\documentclass[11pt,a4paper,english,reqno]{amsart}
\usepackage{amsmath,amssymb,amsfonts,epsfig,mathrsfs}
\usepackage[T1]{fontenc}

\usepackage{color}
\usepackage{array}
\usepackage{amsthm}
\usepackage{amstext}
\usepackage{graphicx}
\usepackage{setspace}
\usepackage[margin=2.5cm]{geometry}
\usepackage{bbm}
\usepackage{color}
\usepackage{enumitem}
\usepackage{undertilde}
\allowdisplaybreaks[4]

\usepackage{pgfplots}

\usepackage{amscd,psfrag}
\usepackage{yhmath}
\usepackage[mathscr]{eucal}

\usepackage{nicefrac}

\usepackage{slashed}

\makeatletter
\pdfpageheight\paperheight
\pdfpagewidth\paperwidth

\usepackage{mathrsfs}

\usepackage{epstopdf}
\usepackage{chngcntr}
\counterwithin{figure}{section}
\usepackage{mathrsfs}

\usepackage{indentfirst}	

\usepackage[normalem]{ulem}
\theoremstyle{plain}

\newtheorem{definition}{Definition}[section]
\newtheorem{theorem}[definition]{Theorem}
\newtheorem*{theorem*}{Theorem}

\newtheorem{remark}[definition]{Remark}

\newtheorem*{remark*}{Remark}
\newtheorem*{sideremark*}{Side Remark}\newtheorem*{mt*}{Main Theorem}

\newtheorem*{claim*}{Claim}
\newtheorem*{q*}{Question}
\newtheorem{lemma}[definition]{Lemma}
\newtheorem{corollary}[definition]{Corollary}
\newtheorem*{corollary*}{Corollary}
\newtheorem*{proposition*}{Proposition}

\newcommand{\R}{\mathbb{R}}

\newcommand{\na}{\nabla}

\newcommand{\dd}{{\rm d}}
\newcommand{\p}{\partial}

\newcommand{\emb}{\hookrightarrow}

\newcommand{\map}{\rightarrow}
\newcommand{\G}{\Gamma}
\newcommand{\M}{{\mathscr{M}}}
\newcommand{\E}{{\mathscr{E}}}
\newcommand{\F}{{\mathscr{F}}}

\newcommand{\dbar}{{\overline{d}}}
\newcommand{\il}{{{\mathbf{i}^{(\ell)}}}}
\newcommand{\ilo}{{{\mathbf{i}^{(\ell-1)}}}}
\newcommand{\pr}{{\mathtt{pr}}}
\newcommand{\prz}{\pr_{\mathtt{H}}}
\newcommand{\prs}{\pr_\star}
\newcommand{\dmue}{\dd \mu^{
\E}}
\newcommand{\cc}{{\mathtt{c}}}

\newcommand{\gge}{{\mathfrak{g}^\E}}

\allowdisplaybreaks[4]

\def\XXint#1#2#3{{\setbox0=\hbox{$#1{#2#3}{\int}$ }
\vcenter{\hbox{$#2#3$ }}\kern-.6\wd0}}

\numberwithin{equation}{section}
\numberwithin{figure}{section}

\title{Towards a Theory of Multi-Parameter Geometrical Variational Problems: Fibre Bundles, Differential Forms and Riemannian Quasiconvexity}

\author{Siran Li}
\address{Siran Li: Department of Mathematics, Rice University, MS 136
P.O. Box 1892, Houston, Texas, 77251-1892, USA\, $\bullet$ \,  Department of Mathematics, McGill University, Burnside Hall, 805 Sherbrooke Street West, Montreal, Quebec, H3A 0B9, Canada.}

\email{\texttt{Siran.Li@rice.edu}}


\date{\today}

\begin{document}

\maketitle

\begin{abstract}
We are concerned with the relaxation and existence theories of a general class of geometrical minimisation problems, with action integrals defined via differential forms over fibre bundles. We find  natural algebraic and analytic  conditions which give rise to a relaxation theory. Moreover, we propose the notion of ``Riemannian quasiconvexity'' for cost functions whose variables are differential forms on Riemannian manifolds, which extends the classical quasiconvexity condition in the Euclidean settings. The existence of minimisers under the Riemannian quasiconvexity condition has been established. This work may serve as a tentative generalisation of the framework developed in the recent paper \cite{dg} by Dacorogna--Gangbo. 
\end{abstract}

\section{Introduction}\label{sec: intro}

In an interesting recent paper \cite{dg}, Dacorogna and Gangbo study the relaxation theory of a family of dynamic variational problems with action integrals defined on paths of differential forms over contractible domains in Euclidean spaces. Such paths are parametrised by a time variable $t \in [0,1]$, thus rendering the dynamic features of the problem in consideration. The existence theory for minimisers has also been established, under a variant of the classical quasiconvexity condition {\it \`{a} la} Morrey \cite{m',m}.

The study in \cite{dg} has been motivated by problems in physics and engineering related to the transport phenomena. Primary examples include the Maxwell equations subject to the ideal Ohm's law ({\it cf}. ``Model Example'' in \S 1, \cite{dg}) and the continuity equation with a generalised kinetic energy ({\it cf}. \S 3.6 {\it op.\,cit.}). The work \cite{dg} completes the programme on relaxation theory in a series of papers; including \cite{bds, cdk, dg'} and many references cited therein. We also refer to the fundamental work \cite{b} by Ball for more about the physical and engineering applications.

The aim of our paper is to discuss, in tentative manners, possible directions in which the relaxation theory programme of Dacorogna--Gangbo  may be further extended. Our attempts embrace of two main features:
\begin{enumerate}
\item
Nontrivial geometry of the spatial domain;
\item
Multi-parameter setting for the variational problem.
\end{enumerate}
$(1)$ means that we consider action integrals over domains on general Riemannian manifolds, not just on Euclidean spaces. The meaning of $(2)$ is that, loosely speaking, we allow the relevant variational problems to have several ``time'' variables $(t_1, \ldots, t_k)$ instead of one single $t$.

Motivated by the above considerations, we propose a model of multi-parameter variational problems over fibre bundles. A fibre bundle is a geometric object that splits into two directions: we view the horizontal direction --- which is a Riemannian manifold termed as the {\em base manifold}, thus of nontrivial,  curved geometry --- as the {\em spatial direction}, and the vertical direction as a $k$-dimensional {\em parameter space}. For the simplest yet fundamental example, the space--time $\R^1 \times \R^3$ is a fibre bundle of spatial domain $\R^3$ and fibre $\R^1$ ($k=1$). In general, the fibres are glued together in a consistent geometrical manner to form the bundle. Therefore, the point $(1)$ in the previous paragraph is modelled by the Riemannian geometry of the base manifolds, and $(2)$ is manifested by the multi-dimensionality of the fibres.

Our investigations in this paper focus on geometrical variational problems {\em in the absence of topological obstructions}. In particular, in the main existence Theorem \ref{thm: existence}, we consider the variational problems over contractible subdomains of fibre bundles. The triviality of topology is, in fact, the working assumption in \cite{dg} and preceding works \cite{bds, d, dg'} etc.. Let us also bring to attention the work \cite{w} by Wang, which is an earlier study of a minimisation problem over fibre bundles, aiming at generalising the analytic theories of harmonic maps to bundles.


Throughout this paper, $\M$ is an $n$-dimensional Riemannian manifold, $\pi: \E \map \M$ is a fibre bundle with typical fibre $\F$ being an $k$-dimensional manifold. (The relevant geometrical backgrounds will be discussed in \S \ref{sec: prelim}.) For given integers $\ell \in \{1,2,\ldots, n\}$ and $I \in \mathbb{N}$, we shall consider a {\em cost function}:
\begin{equation*}
\cc: \Omega^\ell(\M) \times \Omega^{a_1}(\M) \times \ldots \times \Omega^{a_I}(\M) \longrightarrow ]-\infty, \infty],
\end{equation*}
with integers $a_i<\ell$ for each $i \in \{1,2,\ldots, I\}$; $\Omega^\bullet(\M)$ denotes the space of differential forms of degree $\bullet$ over $\M$. The goal of our {\em variational problem} is to minimise the total cost (known as the {\em action integral} in the  calculus of variations):
\begin{equation}\label{A problem, intro}
\mathscr{A}(f;g_1,\ldots,g_I) := \int_{\mathcal{O}} \cc(f, g_1, \ldots, g_I)\,\dmue,
\end{equation}
where $f: \mathcal{O} \subset \E \map \Omega^\ell(\M)$ and $g_i:  \mathcal{O} \map \Omega^{a_i}(\M)$ for each $i$. Here we integrate with respect to $\dmue$, the volume measure on $\E$. Also, as explained in the preceding paragraph, $\mathcal{O}$ is taken to be a contractible subdomain of the bundle $\E$, {\it i.e.}, it is topologically equivalent to a single point. For example, each convex domain (and more generally, each star-shaped domain) is contractible. We emphasise that the variables $\{f; g_1, \ldots, g_I\}$ are taken to be functions over the bundle $\E$, not just on the base manifold $\M$; in other words, they depend on the $k$ parameters over the fibre $\F$.

The main result of this paper, Theorem \ref{thm: existence},  is to specify a natural sufficient condition for the existence of minimisers of the action in \eqref{A problem, intro}. This will be done subject to suitable {\em analytic conditions} concerning the regularity of the variables $f; g_1, \ldots, g_I$ and their boundary conditions on $\p\mathcal{O}$. Our condition, termed as ``{\em Riemannian quasiconvexity}'' in Definition \ref{def: quasiconvex, mfd}, is a generalised version of {\em quasiconvexity} in Morrey's classical theory (\cite{m, m'}). The only difference is that, in order to account for the nontrivial geometry of the fibre bundle, we introduce a volume growth factor ({\it cf}. Definition \ref{def: growth factor}). This factor reduces to unity in the case of Euclidean geometry, hence thereof the notions of Riemannian quasiconvexity and quasiconvexity coincide.

However, Morrey's theory of quasiconvexity is specifically designed for variational problems of the following form (for simplicity, let us now focus on the Euclidean case):
\begin{equation}\label{A' problem, intro}
\mathscr{A}'(u) := \int_{U} \cc'(\na u)\,\dd x,
\end{equation}
where $U \subset \R^N$ is a Euclidean domain and $u: U \map \R^M$. The salient feature is that the cost function $\cc'$ is a function of a gradient map. Thus, we need to recast the variational problem in \eqref{A problem, intro} into the form of \eqref{A' problem, intro}. 

Such a ``translation'' \eqref{A problem, intro} $\rightsquigarrow$ \eqref{A' problem, intro} is known as {\em relaxation}.  Dacorogna--Gangbo developed in \cite{dg} a relaxation theory for actions defined on trivial bundles over Euclidean domains with typical fibre $\R$. In our case, roughly speaking, a relaxation theory amounts to finding a differential form on the bundle which encompasses all the information of the differential forms $\{f; g_1, \ldots, g_I\}$ on the base manifold $\M$, as well as a suitable ``gauged cost function'' $\cc' \equiv \cc_{\rm gauge}$ obtained from $\cc$.

In this paper, we single out an {\em algebraic condition} in terms of $\{f; g_1, \ldots, g_I\}$ that ensures the existence of a relaxation theory. Heuristically, our condition says that $\{f; g_1, \ldots, g_I\}$ arise from the horizontal projection ({\it i.e.}, the morphism between  differential forms naturally induced by the bundle map $\pi: \E \map \M$) of one single exact differential $\ell$-form on $\E$. In other words, all the information of the data $\{f; g_1, \ldots, g_I\}$ for our variational problem \eqref{A problem, intro} can be packaged into a nice object ``upstairs'' on $\E$. In Definition \ref{def: horizontal shadow}, we give this condition a suggestive name: $\{f; g_1, \ldots, g_I\}$ a said to be a ``{\em horizontal shadow}'' of the exact form on $\E$. Under such condition a relaxation theory can be established for the minimisation problem of the action \eqref{A problem, intro}, which, together with the earlier discussions, admits a minimiser provided that the Riemannian quasiconvexity condition is also satisfied.

\noindent
{\bf Organisation for the remaining sections:} In $\S \ref{sec: prelim}$ we present some background materials on differential geometry. In $\S \ref{sec: relaxation}$ we establish the relaxation theory, with focuses on the algebraic aspects. The existence theory under the Riemannian quasiconvexity condition is established in $\S \ref{sec: exist}$. In particular, the main result of the paper is Theorem \ref{thm: existence}. Finally, in $\S \ref{sec: final}$ we conclude with some questions for further studies.

\noindent
{\bf Remark on notations:}  Our notations in this paper are mostly consistent with those in \cite{dg}. In particular, the variables $f, g_1, \ldots, g_I$ denote differential forms. The only major difference is that we use the geometric notation $\Omega^\bullet(\M)$ to denote the space of differential forms --- {\it i.e.}, the space of sections of alternating algebras over the cotangent bundle $T^*\M$ ---  instead of the geometric measure theoretic notation $\Lambda^\bullet(\R^n)$ (for the special case $\M=\R^n$) in \cite{dg}. Moreover, in this paper $g$ always denotes a differential-form-valued function, while $\mathfrak{g}$, $\mathfrak{g}^\E$ denote Riemannian metrics on $\M$ and $\E$, respectively.

\section{Preliminaries}
\label{sec: prelim}

In this section we collect some background materials on differential geometry, with emphases on vector bundles and differential forms. We refer to Part $1$ of the classical text \cite{bt} by Bott--Tu for preliminaries on topology and geometry, and to \cite{h} by Hebey for materials on analysis. A comprehensive treatment of multilinear algebra can also be found in Chapter $1$ of Federer \cite{f}.

Throughout this paper, $\M$ is an $n$-dimensional Riemannian manifold without boundary. A {\em fibre bundle} $\pi: \E \map \M$  consists of the {\em total space} $\E$, the {\em base manifold} $\M$ and the {\em typical fibre} $\F$. Here, $\E$ is an $(n+k)$-dimensional  manifold with Riemannian metric $\gge$, and $\F$ is a $k$-dimensional manifold.  The {\em bundle projection} is a submersion satisfying the local trivialisation requirement: given any point $z\in \M$, there exists an open chart $V_\alpha \subset \M$ in the atlas of $\M$ which contains $z$, such that there is a diffeomorphism $\Phi: \pi^{-1}(V_\alpha) \map V_\alpha \times F$. We call $T\M$ the {\em horizontal} direction and $T\F$ the {\em vertical} direction of the bundle. 

The simplest example of a fibre bundle is the trivial bundle $\E = \M \times \F$. In this case the bundle is globally trivialised, {\it i.e.}, in the above one may choose the chart $V_\alpha$ to be $\M$ for any $z \in \M$. When $\F=\R^k$, $\pi:\E\map\M$ is known as a {\em vector bundle}. By an abuse of notations, we also refer to the total space $\E$ as the bundle. A {\em section} of $\E$ is a smooth map $\beta: \M \map \E$ satisfying $\pi \circ \beta = {\bf Id}_\M$. The space of sections of $\E$ is denoted as $\G(\E)$. 

For simplicity, we assume in this paper that all the Riemannian metrics are at least bounded in $C^2$, and that all the diffeomorphisms  are $C^\infty$.

Given a Riemannian manifold $\M$, we can define on it the space of differential forms. In brief, let $\bigwedge^\ell T^*\M$ be the $\ell^{\text{th}}$-grading of the exterior algebra over $\M$, {\it i.e.}, the vector space of alternating $\ell$-forms over the cotangent bundle $T^*\M$. This is a vector bundle over $\M$. We define
\begin{equation*}
\Omega^\ell (\M) := \G\Big(\bigwedge^\ell T^*\M\Big),
\end{equation*}
and an element of it is known as a {\em differential $\ell$-form}. It is clear that $\Omega^0(\M) = C^\infty(\M;\R)$ and that $\Omega^\ell(\M) = \{0\}$ for $\ell > n =\dim\M$. Using the canonical duality between $T^*\M$ and $T\M$, we see that a $1$-form is canonically due to a vectorfield, {\it i.e.}, an element of $\G(T\M)$. Similarly, we can define differential forms on $\E$ and $\F$.

In a locally trivialised chart of $\E$, there is a coordinate system $\{dx^1, \ldots, dx^n; dx^{n+1}, \ldots, dx^{n+k}\}$ such that $\{dx^1, \ldots, dx^n\}$ spans $\G(T^*\M)$ and that $\{dx^{n+1}, \ldots, dx^{n+k}\}$ spans $\G(T^*\F)$. Thus, any differential $\ell$-form $\xi$ on $\E$ locally takes the form
\begin{equation*}
\xi = \sum_{1 \leq i_1 < \ldots < i_\ell \leq n+k} \xi_{i_1\ldots i_{\ell}} dx^1 \wedge \ldots \wedge dx^{n+k},
\end{equation*}
where $\xi_{i_1\ldots i_{\ell}} \in \R$ is a smooth function on the chart. Each  ordered $\ell$-tuple of indices $(i_1, \ldots, i_\ell)$ in the summation is called a {\em multi-index} of valency $\ell$. By the definition of differential forms, the coefficient $\xi_{i_1\ldots i_{\ell}}$ changes its sign each time a pair of indices gets interchanged; thus we agreed on labelling the indices in the ascending order. 

On each differential manifold $\M$, there is a map on the graded algebra of differential forms, known as the {\em differential}: For each $\ell\in\mathbb{N}$, we have $d: \Omega^\ell(\M) \map \Omega^{\ell
+1}(\M)$ which is multi-linear, satisfies the Leibniz rule, and verifies $d\circ d =0$. For $\varphi \in \Omega^0(\M) = C^\infty(\M;\R)$, the $1$-form $d\varphi$ is dual to the vectorfield $\na \varphi$, which is the gradient of $\varphi$. Throughout this paper, the differential on $\E$ is denoted by $\dbar$. 

Let $\M$ and $\M'$ be two differential manifolds. A diffeomorphism $\Phi: \M \map \M'$ is a smooth map with smooth inverse. Given a differential $\ell$-form $\gamma$ on $\M'$, we can {\em pullback} the form to get $\Phi^\#\gamma \in \Omega^\ell(\M)$. Similarly, we can pushforward form $\M$ to $\M'$ by setting $\Phi_\# := (\Phi^{-1})^\#$. 

Over a Riemannian manifold we may define the {\em Sobolev spaces} $W^{1,s}$ of differential forms. This can be done intrinsically, {\it i.e.}, only using the Riemannian structure of the manifold. In what follows we shall focus on the range $1<s<\infty$. For $V \subset \M$ and $\E'$ a fibre bundle over $\M$, we use the notation $W^{1,s}(V,\E')$ to denote the space of $W^{1,s}$-sections of $\E'$ defined over $V$.

One can integrate differential $0$-forms or $n$-forms on $n$-dimensional manifolds over the {\em Riemannian volume measure} induced from the metric. In this paper, $\dmue$ denotes the volume measure on $\E$, and $\dd \mathcal{L}^N$ denotes the $N$-dimensional Lebesgue measure on Euclidean spaces.

Finally, let us introduce the {\em comass norm} on differential forms on $\M$. Let $\mathfrak{g}$ the Riemannian metric on $\M$. A section $v \in \G(\bigwedge^\ell T\M)$ is said to be {\em simple} if it equals a field of alternating products of $\ell$ elements in $T\M$. Naturally one may extend the metric $\mathfrak{g}$ to a field of inner products between $\bigwedge^\ell T\M$ and $\bigwedge^\ell T^*\M$, and hence to their sections. Then, the {\em comass norm} of $\psi \in \Omega^\ell(\M)$ is defined by
\begin{equation*}
\|\psi\| := \sup\Bigg\{ \mathfrak{g}(\psi, \chi): \, \chi \in \G\Big( \bigwedge^\ell T\M\Big), \, \chi\text{ is simple, and } \mathfrak{g}(\chi, \chi) \leq 1 \Bigg\};
\end{equation*}
see Federer \cite{f}, 1.8.1. The term $\mathfrak{g}(\chi, \chi)$ is understood with the obvious duality.

\section{Relaxation}\label{sec: relaxation}

In this section, we develop the relaxation theory for the minimisation problem of the action integral \eqref{A problem, intro} in the Introduction.

\subsection{Algebraic Condition}
In this first step, we exhibit an algebraic condition that ensures the reduction of \eqref{A problem, intro} to \eqref{A' problem, intro} (Theorem \ref{thm: main}). At this moment we do not impose any further regularity requirement; everything is assumed to be smooth. Nevertheless, the constructions in this subsection can be easily generalised to the case $\{f; g_1, \ldots, g_I\} \in L^s$ for $s \in ]1,\infty[$, via standard approximation techniques.

\begin{definition}\label{def: l-complementary}
Let $\ell \in \mathbb{N}$ and let $g \in \Omega^{j}(\M)$, $0\leq j \leq \ell$, be a differential $j$-form on $\M$. Any differential $(\ell-j)$-form defined on the fibre $\F$ is said to be \underline{$\ell$-complementary} to $g$. 
\end{definition}
The above definition is motivated by a na\"{i}ve observation: if $\vartheta$ is $\ell$-complementary to $g$, then $g \wedge \vartheta$ is an $\ell$-form on the total space $\E$, whose horizontal part is precisely $g$. In local coordinates, if $\{d\eta^1, \ldots, d\eta^{k}\}$ is a coframe for $T^*\F$, then we may express $\vartheta = d\eta^{i_1}\wedge\ldots\wedge d\eta^{i_{\ell-j}}$ with $1 \leq i_1 \leq \ldots \leq i_{\ell-j} \leq k$. In what follows we are only interested in the case when $\vartheta$ is closed. 

\begin{definition}\label{def: horizontal shadow}
Let $\mathscr{G}:=\{g_1, \ldots, g_I\}$ be a finite collection of differential forms of degree $\leq \ell$ on $\M$ and let $\varpi \in \Omega^{\ell}(\E)$. We say that $\mathscr{G}$ consists of \underline{horizontal shadows of $\varpi$} if and only if the following holds: Suppose that in a local coordinate system $\{dx^1, \ldots dx^{n+k}\}$ for $\E$ one has 
\begin{equation*}
\varpi = \sum_{1 \leq i_1 < \ldots < i_\ell \leq n+k} \varpi_{i_1 \ldots i_\ell} dx^{i_1} \wedge \ldots \wedge dx^{i_\ell}.
\end{equation*}
Then, for each $g_j \in \mathscr{G}$, there are a multi-index $\il=(i_1, \ldots, i_\ell)$ and an index $\star\in \{1,2,\ldots,\ell\}$ depending on $j$ such that 
\begin{equation*}
g_j = \varpi_{i_1 \ldots i_\ell} dx^{i_1}\wedge\ldots\wedge dx^{i_\star}.
\end{equation*}
\end{definition}

In other words, by a ``horizontal shadow'' we mean the horizontal projection of some component of $\varpi$. We emphasise that the members of $\mathscr{G}$ may have values depending on the fibre $\F$, though they are genuine differential forms over the base manifold $\M$.

Our major concern is to seek criteria for $\{g_1, \ldots, g_I\} \subset \bigcup_{j \leq \ell}\Omega^j(\M)$ to be {\em simultaneous} horizontal shadows of one exact differential $\ell$-form on $\E$. For this purpose, we  analyse the linear combinations of $g_j \wedge \vartheta^j$, where $\vartheta^j$ are $\ell$-complementary to $g_j$ as in Definition \ref{def: l-complementary}.

\begin{theorem}\label{thm: main}
Let $\pi: \E^{n+k} \map \M^{n}$ be a fibre bundle with $\F$ denoting the typical fibre.  Let $\dbar$ be the exterior differential on $\E$. Consider the following data:
\begin{itemize}
\item
A differential form $f \in \Omega^{\ell}(\M)$;
\item
A collection of differential forms $\{g_1, \ldots, g_I\}$ on $\M$, each of degree no more than $(\ell-1)$.
\end{itemize}
\begin{enumerate}
\item
Suppose that for some closed forms $\vartheta^1, \ldots, \vartheta^I$ which are $\ell$-complementary to $g_1, \ldots, g_I$, respectively, there holds
\begin{equation}\label{closedness assumption}
\dbar\Big(f + \sum_{i=1}^I  g_i \wedge \vartheta^i\Big) = 0.
\end{equation}
Then one can find a differential form $\xi \in \Omega^{\ell-1}(\E)$ on the total space and a closed form $h_0 \in \Omega^\ell(\E)$, such that the data $\{f; g_1, \ldots, g_I\}$ are horizontal shadows of $\dbar\xi + h_0$. 
\item
Conversely, for each $\xi \in \Omega^{\ell-1}(\E)$ we can find $f\in \Omega^{\ell}(\M)$ and a collection of forms $\{g_1, \ldots, g_I\}$ on $\M$, each $g_i$ of degree no more than $(\ell-1)$, such that \eqref{closedness assumption} holds for some $\vartheta^1, \ldots, \vartheta^I$ that are $\ell$-complementary to $g_1, \ldots, g_I$, respectively, and that $\{f; g_1, \ldots, g_I\}$ are horizontal shadows of $\dbar\xi$.
\end{enumerate}
\end{theorem}

\begin{remark}
In the setting of $(1)$ above,  let us write 
\begin{equation}\label{prh-notation}
\prz\big(\dbar\xi+h_0\big) = (f; g_1, \ldots, g_I).
\end{equation}
Notice that $\prz$ is a well-defined map from $\Omega^\ell(\E)$ to $\bigcup_{j=1}^\ell \Omega^j(\M)$, up to the permutations of $g_i$'s; it is naturally induced by the bundle projection $\pi: \E \map \M$.
\end{remark} 

\begin{proof}
	The first statement is simple: thanks to the condition \eqref{closedness assumption}, $f+\sum_{i}g_i \wedge \vartheta^i$ is a  representative of $H^{\ell}(\E;\R)$ = the $\ell^{\text{th}}$-de Rham cohomology group of $\E$. 
Thus, there exists a closed $\ell$-form $h_0$ such that $$f+\sum_{i}g_i \wedge \vartheta^i - h_0 = \dbar\xi$$ for some $\xi \in \Omega^{\ell-1}(\E)$. The statement then follows from Definition \ref{def: horizontal shadow} of the horizontal shadows. 
 	
 	To prove the second statement, it suffices to assume that $\xi$ is not purely vertical, {\it i.e.}, $\xi \in \Omega^{\ell-1}(\E) \setminus [\Omega^0(\M)\otimes \Omega^{\ell-1}(\F)]$, since in this case one may take $f \equiv g_i \equiv 1$. 	Next, let us express $\xi$ in the local coordinates $\{dx^1, \ldots, dx^{n+k}\}$, whereof the first $n$-coordinates form a coframe for $T^*\M$, and the last $k$-coordinates form a coframe for $T^*\F$. That is, 
 	\begin{align*}
 	\xi = \sum_{1\leq i_1 < \ldots < i_{\ell-1}\leq n+k} \, \xi_{i_1\ldots i_{\ell-1}}dx^{i_1} \wedge \ldots \wedge dx^{i_{\ell-1}}
 	\end{align*}
 hence
 \begin{align*}
 \dbar\xi = \sum_{j=1}^{n+k} \,\,\sum_{1\leq i_1 < \ldots < i_{\ell-1}\leq n+k} \, \frac{\p\xi_{i_1\ldots i_{\ell-1}}}{\p x^j}dx^j \wedge dx^{i_1} \wedge \ldots \wedge dx^{i_{\ell-1}}.
\end{align*}
Let us denote by $\star \in \{1, 2, \ldots, \ell-1\}$ the largest integer such that $i_\star \leq n$. One may express
\begin{align*}
\dbar\xi = \bigg( \sum_{j=1}^n + \sum_{j=n+1}^{n+k} \bigg) \Bigg\{ \sum_{\star=1}^{\ell-1} \bigg( \sum_{1 \leq i_1 < \ldots < i_\star \leq n}+ \sum_{n+1 \leq i_{\star+1} < \ldots < i_{\ell-1} \leq n+k} \bigg)  \frac{\p\xi_{i_1\ldots i_{\ell-1}}}{\p x^j}dx^j \wedge dx^{i_1} \wedge \ldots \wedge dx^{i_{\ell-1}} \Bigg\}.
\end{align*}
Here and throughout, the convention is that $\sum_{n+1 \leq i_{\star+1} \leq \ldots \leq i_{\ell-1} \leq n+k} \{\bullet\} \equiv 0$ for $\star = \ell-1$. We start the summation from $\star=1$ since $\xi$ is not purely vertical. 

Now, define the horizontal form $f \in \Omega^{\ell}(\M)$ by
	\begin{equation*}
	f := \sum_{j=1}^n  \,\,\sum_{1 \leq i_1 < \ldots < i_{\ell-1} \leq n} \frac{\p\xi_{i_1\ldots i_{\ell-1}}}{\p x^j}dx^j \wedge dx^{i_1} \wedge \ldots \wedge dx^{i_{\ell-1}}.
	\end{equation*}
Moreover, for each fixed multi-index $\ilo=(i_1, \ldots, i_{\ell-1})$ (hence with $\star \in \{1,2,\ldots,\ell-2\}$ fixed too), define 
\begin{equation*}
g_{\{\ilo, j\}} := \begin{cases}
\frac{\p\xi_{i_1\ldots i_{\ell-1}}}{\p x^j}dx^j \wedge dx^{i_1} \wedge \ldots \wedge dx^{i_{\star}} \qquad \text{ if } 1 \leq j \leq n,\\
\frac{\p\xi_{i_1\ldots i_{\ell-1}}}{\p x^j} dx^{i_1} \wedge \ldots \wedge dx^{i_{\star}} \qquad \text{ if }  n+1 \leq j \leq n+k.
\end{cases}
\end{equation*}
These are differential forms on $\M$ of degree no more than $(\ell-1)$. Let us also put
\begin{equation*}
\vartheta^{\{\ilo, j\}} := \begin{cases} 
dx^{i_{\star+1}} \wedge \ldots \wedge dx^{i_{\ell-1}} \qquad \text{ if }  1 \leq j \leq n,\\
(-1)^\star dx^j \wedge dx^{i_{\star+1}} \wedge \ldots \wedge dx^{i_{\ell-1}}  \qquad \text{ if }  n+1 \leq j \leq n+k\\
\end{cases}
\end{equation*}
By a standard glueing argument in differential geometry, all the differential forms $f$, $g_{\{\ilo, j\}}$ and $\vartheta^{\{\ilo, j\}}$ introduced above can be defined globally on $\M$ or $\F$. In addition, they can be viewed as differential forms on the total space $\E$ via the natural inclusion maps.

By our definition of the index $\star$, the closed form $\vartheta^{\{\ilo,j\}}$ is $\ell$-complimentary to $g_{\{\ilo, j\}}$. It follows that
\begin{align}\label{dbar-xi, 1}
\dbar \xi = f + \sum_{j=1}^{n+k} \,\,\sum_{\star=1}^{\ell-2} \bigg( \sum_{1 \leq i_1 < \ldots < i_\star \leq n}+ \sum_{n+1 \leq i_{\star+1} < \ldots < i_{\ell-1} \leq n+k} \bigg) g_{\{\ilo, j\}} \wedge\vartheta^{\{\ilo, j\}}.
\end{align}
We can now relabel $g_i \equiv g_{\{\ilo, j\}}$ and $\vartheta^i \equiv \vartheta^{\{\ilo,j\}}$ with the new indices
\begin{align*}
i \in \mathcal{I} &:=\Big\{ \big(j, \star, \ilo\big):  1\leq j \leq n+k,\nonumber\\
&\qquad\qquad 1 \leq \star \leq  \ell-2, \, 1\leq i_1 < \ldots < i_\star \leq n \text{ and } n+1 \leq i_{\star+1} < \ldots < i_{\ell-1} \leq n+k \Big\}.
\end{align*}
The indexing set $\mathcal{I}$ is finite; in fact, 
\begin{equation*}
{\rm card}(\mathcal{I}) \leq (n+k) \sum_{\star=1}^{\ell-2} {{n}\choose{\star}} \cdot  {{k}\choose{\ell-\star}} \leq C(n,k)<\infty.
\end{equation*}

In summary, Eq.\,\eqref{dbar-xi, 1} can be recast into a finite sum:
\begin{equation*}
\dbar\xi = f + \sum_{i \in \mathcal{I}} g_i \wedge \vartheta^i.
\end{equation*}
It is automatically closed on $\E$. The proof is now complete.   \end{proof}

\begin{remark}
If $H^{\ell}(\E;\R)=\{0\}$, {\it i.e.}, the $\ell^{\text{th}}$-de Rham cohomology group of $\E$ is trivial, then $h_0$ in Theorem \ref{thm: main}(1) can be taken as zero. This is clear from the proof. 
\end{remark}

\subsection{Analytic Conditions}
 
The statement and proof of Theorem \ref{thm: main} above are purely based on multilinear algebraic  computations. In the actual problems arising from physical and engineering applications, one often encounters further {\em regularity} and {\em boundary} conditions, which are analytic in nature. These shall be taken into considerations in this subsection. 

From now on, let $\mathcal{O}$ be a smooth, bounded subdomain of $\E$ and let $I$ be a finite positive integer. Without loss of generality, it will be assumed that $\pi|\mathcal{O}$ still surjects onto $\M$. We consider cost functions of the form 
\begin{equation}\label{cost}
\cc:\Omega^{\ell}(\M) \times \Big[\bigcup_{j=1}^\ell\Omega^j(\M)\Big]^I \longrightarrow ]-\infty,\infty],
\end{equation}
as well as the corresponding action integral: 
\begin{equation}\label{action}
\mathscr{A}\big(f; g_1,\ldots,g_I\big) := \int_{\mathcal{O}} \cc\big(f; g_1, \ldots, g_I\big)\,\dmue.
\end{equation}
The argument $f$ is parametrised by the fibre $\F$, namely $f=f(\sigma, x) \in \Omega^{\ell}(\M)$ with $x \in \M$ and $\sigma \in \F$; similarly for $g_1, \ldots, g_I$. That is, for suitable $\bullet \in \{0,1,\ldots,\ell\}$, we have
\begin{equation*}
f; g_1, \ldots, g_l: \mathcal{O} \subset \E \longrightarrow \Omega^\bullet(\M)
\end{equation*}
as functions.

Fix $1<s<\infty$. 
By a {\em gauge form} (``gauge'' in brief) we mean $$\widetilde{\xi}\in W^{1,s}\Big(\mathcal{O}; \bigwedge^{\ell}T^*\E\Big),$$ {\it i.e.}, a differential $\ell$-form on $\mathcal{O}$ with designated Sobolev regularity. As shall be clear below, gauges are instrumental to the specification of boundary conditions.

In the case $\M=\R^3, \F = \R^1$ and $\E=\M \times \F = \R^4$, a canonical choice for the electro-magnto-dynamical problems is $\widetilde{\xi}:=(\p_t\varphi, \na_x \wedge A)$, where $\varphi, A$ are the scalar and vector electromagnetic potentials, respectively. This is known as the {\em Lorenz gauge}. 
\begin{definition}\label{def: Ps tilde-xi}
The \underline{admissible class  of $L^s$-regularity subject to the gauge $\widetilde{\xi}$} is 
\begin{align}
\wp^s(\widetilde{\xi}) &:= \Bigg\{ (f;g_1, \ldots, g_I) \in L^s\bigg(\mathcal{O}; \bigwedge^{\ell}T^*\M\times \Big[\bigcup_{j=1}^\ell\bigwedge^j T^*\M\Big]^I\bigg): \nonumber\\
&\qquad \qquad  \text{ There are $\ell$-complementary forms } \vartheta^1,\ldots,\vartheta^I \text{ such that } f+ \sum_{i=1}^I g_i \wedge \vartheta^i \text{ is closed, }\nonumber\\
&\qquad\qquad \text{ and that } \Big(f+ \sum_{i=1}^I g_i \wedge \vartheta^i + \dbar\widetilde{\xi}\Big)\big|\p\mathcal{O} \in \G\big(T(\p\mathcal{O})\big)  \Bigg\}.
\end{align}
\end{definition}

The last condition means that $$\gge\bigg(\nu,\Big(f+ \sum_{i=1}^I g_i \wedge \vartheta^i + \dbar\widetilde{\xi}\Big)\bigg)= 0\qquad \text{ on } \p\mathcal{O},$$ where $\nu$ is the outward unit normal vectorfield along $\p\mathcal{O}$ and 
$\gge$ is the bundle metric on $\E$. By our assumptions in $\S \ref{sec: prelim}$, $\nu$ is a smooth vectorfield. Here and hereafter, we shall always identify $\nu$ with the $1$-form obtained via the canonical duality $T\E \cong T^*\E$. Although the trace may fail to be well-defined for differential forms merely of $L^s$-regularity, the above tangency condition nevertheless makes sense. The closedness of differential forms in the other condition will also be understood in the weak ({\it i.e.}, distributional) sense.
  
Taking into account the relevant regularity and boundary conditions in Theorem \ref{thm: main}, we can easily obtain the following relaxation theorem. It extends Proposition 2.7 in \cite{dg} to the setting of general fibre bundles:
\begin{theorem}[$(\wp) = (\wp_{\rm gauge})$]\label{thm: equivalence of variational problems}
Let $\pi: \E^{n+k} \map \M^{n}$ be the fibre bundle with $\F$ denoting the typical fibre.  Let $\dbar$ be the exterior differential on $\E$. Assume that $\mathcal{O} \Subset \E$ is a smooth, contractible, bounded open subdomain. The following variational problems are equivalent:
\begin{equation*}
\inf\Bigg\{ \int_{\mathcal{O}} \cc\Big( f(x,\sigma); g_1(x,\sigma), \ldots, g_I(x,\sigma) \Big) \,\dmue(x,\sigma) :\, (f; g_1, \ldots, g_I) \in \wp^s (\widetilde{\xi}) \Bigg\} \tag{$\wp$} 
\end{equation*}
and 
\begin{equation*}
\inf\Bigg\{ \int_{\mathcal{O}} \cc \circ \prz \big(\dbar\xi(z)\big) \,\dmue(z) :\, \xi - \widetilde{\xi}\in W^{1,s}_0\big( \mathcal{O}; \bigwedge^{\ell}T^*\E\big) \Bigg\}.\tag{$\wp_{\rm gauge}$}
\end{equation*}
Here and throughout, $W^{1,s}_0$ denotes the Sobolev space of $W^{1,s}$-regularity with trace zero. 
\end{theorem}

Employing the notations in \cite{dg}, let us also write
\begin{equation*}
\begin{cases}
\cc_{\rm gauge} := \cc \circ \prz,\\
\wp^s_{\rm gauge}(\widetilde{\xi}) := \widetilde{\xi} + W^{1,s}_0\big( \mathcal{O}; \bigwedge^{\ell}T^*\E\big),
\end{cases}
\end{equation*}
which shall be understood as the ``gauged'' versions of the cost function and the admissible class. Therefore, we have succeeded in developing a relaxation theory, {\it i.e.}, translating the original problem $(\wp)$ to the gauged problem $(\wp_{\rm gauge})$. The form of the latter problem is classical to the study of  calculus of variations; see \cite{m} by Morrey:

\begin{proof}
It essentially follows from Theorem \ref{thm: main}; we only need to check the relevant boundary and regularity conditions. 

First let us take $(f;g_1, \ldots, g_I) \in \wp^s(\widetilde{\xi})$. By Theorem \ref{thm: main} (1) and the ensuing remark, there is $\xi \in \Omega^{\ell-1}(\E)$ such that $\prz(\dbar\xi)=(f;g_1, \ldots, g_I)$ and that $\dbar \xi = f+ \sum_{i=1}^I g_i \wedge \vartheta^i$. We want to further impose $\dbar(\xi + \widetilde{\xi}) |\p\mathcal{O} \in \G(T(\p\mathcal{O}))$. This can be achieved by solving the following boundary value problem with an application of the Stokes' theorem; see, {\it e.g.}, \cite{s} by Schwarz: 
\begin{equation*}
\begin{cases}
-\dbar\xi = f+ \sum_{i=1}^I g_i \wedge \vartheta^i\qquad \text{ in } \mathcal{O},\\
\xi = -\widetilde{\xi}\qquad \text{ on } \p\mathcal{O}.
\end{cases}
\end{equation*}

Conversely, with $\xi \in \Omega^{\ell-1}(\E)$ given, from Theorem \ref{thm: main} $(2)$ one may find a collection of horizontal shadows $\{f; g_1, \ldots, g_I\}$. In particular, $f+ \sum_{i=1}^I g_i \wedge \vartheta^i =:h$ is closed. We can thus solve  (see \cite{s}) for $\xi \in W^{1,s}(\mathcal{O}; \Omega^{\ell-1}(\E))$ from the system: \begin{equation*}
\begin{cases}
\dbar\xi = -h \qquad \text{ in } \mathcal{O},\\
\xi = -\widetilde{\xi}\qquad \text{ on } \p\mathcal{O}.
\end{cases}
\end{equation*}
Again, by virtue of the Stokes' theorem, it implies that $
(f; g_1, \ldots, g_I) \in \wp^s(\widetilde{\xi})$. The proof is now complete.   \end{proof}

\section{Existence}\label{sec: exist}

We are now concerned with the existence of the minimisers for the variational problem $(\wp)$. 

By Morrey's classical theory \cite{m'}, for a variational problem of the form
\begin{equation}\label{A' problem}
\mathscr{A}'(u) := \int_U F\big(x, \na u(x)\big)\,\dd x
\end{equation}
where $U \subset \R^N$ is a bounded open set, $u: U \map \R^M$ and $F \in C^0(U \times (\R^{M}\otimes \R^{N}))$, the lower semicontinuity of $\mathscr{A}'$ ({\it e.g.}, with respect to the weak $W^{1,s}$-topology for $s \in ]1,\infty[$) is, in the most general case, equivalent to the {\em quasiconvexity} of $F$. Let us also remark that the continuity condition on $F$ can be relaxed to Borel measurability plus local boundedness. We recall:
\begin{definition}\label{def: quasiconvex, Euclidean}
In the setting of \eqref{A' problem}, $(x,p) \mapsto F(x,p)$ is a quasiconvex function if and only if 
\begin{equation*}
\frac{1}{\mathcal{L}^{N}(\mathscr{D})}\int_{\mathscr{D}} F\big(x, p + \na\zeta(y)\big)\,\dd y \geq F(x,p)
\end{equation*}
for each fixed $x \in U$, each fixed $p \in \R^M\otimes \R^N$, each domain $\mathscr{D}\subset U$ and each $\zeta \in W^{1,\infty}_0(\mathscr{D}; \R^M)$.
\end{definition}

One runs into difficulties  when trying to generalise this definition to Riemannian manifolds. Roughly speaking, the condition in Definition \ref{def: quasiconvex, Euclidean} has been found using Euclidean congruences (translations and dilations; see \S 2 of \cite{m'}), which cannot be directly extended to manifolds by simply changing the Lebesgue measure to the Riemannian volume measure. To bypass this obstacle, we introduce a geometrical factor to account for the nonlinearities caused by the manifold geometry:
\begin{definition}\label{def: growth factor}
Let $\pi: \E^{n+k} \map \M^{n}$ be a fibre bundle with $\F$ denoting the typical fibre. Assume that $\mathcal{O} \Subset \E$ is a smooth, contractible, bounded open subdomain.  Denote by $g^\E$ and $\dmue$ the Riemannian metric and the volume form on $\E$, respectively. Let $\mathscr{D}\subset\mathcal{O}$ be a subdomain and let $x_0 \in \mathscr{D}$. The \underline{volume growth factor for $\mathscr{D}$ at $x_0$} is
\begin{equation*}
\mathcal{V}(x_0, \mathscr{D}) := \int_{\mathscr{D}} \frac{\sqrt{\det \,g^\E(x_0)}}{\sqrt{\det \,g^\E(x)}} \,\dmue(x).
\end{equation*}
\end{definition}

For our purpose, let us consider the following variant of Definition \ref{def: quasiconvex, Euclidean}:
\begin{definition}\label{def: quasiconvex, mfd}
Let $\pi, \E,\M,\F$ and $\mathcal{O}$ be as in Definition \ref{def: growth factor}.  Let $\dbar$ be the exterior differential on $\E$. Denote by $\gge$ and $\dmue$ the Riemannian metric and the volume form on $\E$, respectively. A continuous function $F: \mathcal{O} \times \Omega^1(\mathcal{O}) \map \R$ is said to be \underline{Riemannian quasiconvex} if and only if
\begin{equation}
\frac{1}{\mathcal{V}(x_0, \mathscr{D})}\int_{\mathscr{D}} F\big(x, p+\dbar\zeta\big)\,\dmue \geq F\big(x_0, p\big)
\end{equation}
for each fixed $x_0 \in \mathscr{D}$, each fixed $p \in \Omega^1(\mathcal{O})$, each subdomain $\mathscr{D} \subset \mathcal{O}$ and each test function $\zeta \in W^{1,\infty}_0(\mathscr{D})$.  
\end{definition}

As suggested by its name, the notion of Riemannian quasiconvexity on $\E$ depends only on the Riemannian structure of $\E$, not on the fibre bundle structure. Thus, this definition well applies to any Riemannian manifold.

In Definition \ref{def: quasiconvex, mfd}, if $\E$ is the Euclidean space $\R^N$, then $\mathcal{V}(x_0, \mathscr{D}) \equiv \mathcal{L}^N(\mathscr{D})$. Thus, ``Riemannian quasiconvexity'' generalises the classical concept of quasiconvexity over Euclidean spaces, {\it i.e.}, Definition \ref{def: quasiconvex, Euclidean}.

In passing, we also comment that the definition of quasiconvexity has been extended to the setting of differential forms over {\em Euclidean} spaces ({\it cf.} \cite{bds, cdk, dg} and the references cited therein), as well as the related notions of  convexity, polyconvexity and rank-one convexity.

Now we can readily present the following existence theorem. Our proof follows from an adaptation of the direct method for calculus of variations ({\it cf.} {\it e.g.},  Theorem 5.1 in \cite{bds}). The norm $\|\bullet\|$ for (a finite collection of) differential forms is the comass norm; see \S \ref{sec: prelim} above.

\begin{theorem}\label{thm: existence}
Let $\pi: \E^{n+k} \map \M^{n}$ be a fibre bundle with $\F$ denoting the typical fibre. Assume that $\mathcal{O} \Subset \E$ is a smooth, contractible, bounded open subdomain. Take any $s \in ]1,\infty[$ and let $\cc$ be a cost function as in \eqref{cost}:
\begin{equation*}
\cc:\Omega^{\ell}(\M) \times \Big[\bigcup_{j=1}^\ell\Omega^j(\M)\Big]^I \longrightarrow ]-\infty,\infty].
\end{equation*}
Assume that $\cc\circ \prz$ is Riemannian quasiconvex (see Definition \ref{def: quasiconvex, mfd}), and that there are constants $a_1 \in \R$ and $a_2, b_1, b_2 \in ]0,\infty[$ satisfying
\begin{equation}\label{coercive condition}
a_1 + b_1\|(f;g_1, \ldots, g_I)\|^s \leq \cc(f; g_1, \ldots, g_I) \leq a_2 + b_2\|(f;g_1, \ldots, g_I)\|^s
\end{equation}
for each $f: \mathcal{O} \map  \Omega^\ell(\M)$ and $g_1, \ldots, g_I: \mathcal{O} \map  \bigcup_{j=1}^\ell \Omega^j(\M)$. Then the variational problem $(\wp)$ has its minimum attained in the affine Sobolev space $[\widetilde{\xi} + W^{1,s}_0( \mathcal{O}; \bigwedge^{\ell}T^*\E)]$.
\end{theorem}

Recall from \eqref{cost} that $\cc_{\rm gauge} := \cc\circ\prz$ is the gauged cost function. Also, \eqref{coercive condition} is known as a {\em coercive condition}.

\begin{proof}
	In view of Theorem \ref{thm: equivalence of variational problems}, namely that $(\wp) = (\wp_{\rm gauge})$ on $\mathcal{O}$, it suffices to work with the gauged problem $(\wp_{\rm gauge})$. Let $\{\omega_j\}$ be a minimising sequence for $(\wp_{\rm gauge})$, {\it i.e.}, 
\begin{align}\label{min}
\int_{\mathcal{O}} \cc \circ \prz (\dbar\omega_j) \,\dmue \longrightarrow\inf\Bigg\{ \int_{\mathcal{O}} \cc \circ \prz (\dbar\xi) \,\dmue :\, \xi - \widetilde{\xi}\in W^{1,s}_0\big( \mathcal{O}; \bigwedge^{\ell}T^*\E\big) \Bigg\}.
\end{align}
Let us define $\{\xi_j\}$ via the following boundary value problem:
\begin{equation}\label{bvp, intermediate}
\begin{cases}
\dbar \xi_j = \dbar\omega_j\qquad \text{ in } \mathcal{O},\\
\dbar^* \xi_j = 0 \qquad \text{ in } \mathcal{O},\\
\gge(\nu, \xi_j) = \gge(\nu, \omega_j) \equiv -\gge(\nu, \widetilde{\xi})\qquad \text{ on } \p\mathcal{O}. 
\end{cases}
\end{equation}
Recall that the gauge form $\widetilde{\xi}\in W^{1,s}(\mathcal{O}; \bigwedge^{\ell}T^*\E)$, and that the boundary condition is ensured by the construction of $(\wp_{\rm gauge})$. From the classical div-curl estimate (see \cite{cdk, s}) one may infer:
\begin{equation*}
\|\xi_j\|_{W^{1,s}} \leq C_1\Big( \|\dbar \omega_j\|_{L^s} + \|\widetilde{\xi}\|_{W^{1,s}}\Big) \leq C_2. 
\end{equation*}
In the above estimate, $C_2$ is a uniform constant, thanks to the coercivity assumption \eqref{coercive condition} on $\cc$. For this purpose, it is crucial to notice that the definition of $\prz$ ensures the equivalence between \eqref{coercive condition} and the corresponding coercive condition on $\cc\circ\prz$.

After passing to a subsequence (not relabelled), $\{\xi_j\}$ converges weakly in $W^{1,s}(\mathcal{O}; \bigwedge^{\ell}T^*\E)$ to some $\xi$. We modify such $\xi \in W^{1,s}(\mathcal{O}; \bigwedge^{\ell}T^*\E)$ by solving (again, see \cite{cdk}) for $\widehat{\xi}$ via
\begin{equation}\label{hat-xi def}
\begin{cases}
\dbar \widehat{\xi} = \dbar \xi\qquad \text{ in } \mathcal{O},\\
\widehat{\xi} = \widetilde{\xi}\qquad \text{ on } \p\mathcal{O}.
\end{cases}
\end{equation}
In particular, 
\begin{equation*}
\int_{\mathcal{O}} \cc \circ \prz (\dbar \widehat{\xi})\,\dmue  \equiv	\int_{\mathcal{O}} \cc \circ \prz (\dbar \xi)\,\dmue.
\end{equation*}	

Now, let us {\em claim} the following {\em lower semicontinuity property}:
\begin{equation}\label{claim: lsc}
\int_{\mathcal{O}} \cc \circ \prz (\dbar \xi)\,\dmue \leq \liminf_{j \map \infty} 	\int_{\mathcal{O}} \cc \circ \prz (\dbar \xi_j)\,\dmue.
\end{equation}
Assuming it for the moment, we notice that the right-hand side of \eqref{claim: lsc} equals 
\begin{equation*}
\inf\Bigg\{ \int_{\mathcal{O}} \cc \circ \prz (\dbar\xi) \,\dmue :\, \xi - \widetilde{\xi}\in W^{1,s}_0\big( \mathcal{O}; \bigwedge^{\ell}T^*\E\big) \Bigg\},
\end{equation*}
by virtue of \eqref{min} and \eqref{bvp, intermediate}. Thus, in view of \eqref{hat-xi def},  $\widehat{\xi}$ attains the minimum of the gauged problem $(\wp_{\rm gauge})$. This shall complete the proof.

It thus remains to establish the {\em claim} \eqref{claim: lsc}. Here we make crucial use of the Riemannian quasiconvexity condition. Indeed, by the assumption for $\mathcal{O}$, there is a diffeomorphism $$\Phi: (U \subset \R^{N\equiv n+k}) \longrightarrow (\mathcal{O} \subset \E).$$ Let us fix $\mathscr{D}\subset\mathcal{O}$, $x_0 \in \mathscr{D}$, ${p} \in \Omega^1(\mathscr{D})$ and $\zeta \in W^{1,\infty}_0(\mathscr{D})$.  Then, by a change of variables, we can express the Riemannian quasiconvexity of $\cc\circ\prz$ as
\begin{align*}
&\frac{1}{\mathcal{V}(x_0, \mathscr{D})} \int_{\Phi^{-1}(\mathscr{D})}\bigg\{ \Phi^\# [\cc\circ\prz]\Big(\Phi(x_0), \Phi^\#p + \Phi^\#\dbar\zeta(y)\Big)\sqrt{\det\,g^\E(y)}\bigg\}\,\dd \mathcal{L}^N(y) \\
&\qquad\qquad\geq \cc\circ\prz(p).
\end{align*}
In addition, note that
\begin{equation*}
{\mathcal{V}(x_0, \mathscr{D})} = \sqrt{\det\,g^\E(x_0)} \int_{\Phi^{-1}(\mathscr{D})} 1\, \dd \mathcal{L}^{N}, 
\end{equation*}
so the above condition is equivalent to
\begin{align}\label{E quasiconvex, equivalent}
&\frac{1}{\mathcal{L}^N\big(\Phi^{-1}(\mathscr{D})\big)}\int_{\Phi^{-1}(\mathscr{D})}\bigg\{\Phi^\# [\cc\circ\prz]\Big(\Phi(x_0), \Phi^\#p + \Phi^\#\dbar\zeta(y)\Big)\sqrt{\det\,g^\E(y)}\bigg\}\,\dd \mathcal{L}^N(y)\nonumber\\
&\qquad\qquad\qquad \geq \sqrt{\det\,g^\E(x_0)}\,\,\cc\circ\prz(p) \nonumber\\
&\qquad\qquad\qquad \equiv \sqrt{\det\,g^\E(x_0)} \,\, \Phi^\# [\cc\circ\prz](\Phi^\#p). 
\end{align}

To proceed, let us note that $\widehat{p}:=\Phi^\# p$ is an arbitrary $1$-form on $\Phi^{-1}(\mathscr{D})$, as $p\in \Omega^1(\mathscr{D})$ is arbitrary and $\Phi$ is a diffeomorphism. In addition, denoting by $d^{\bf E}$ the exterior differential on the Euclidean space $\R^N$, we have
\begin{equation*}
\Phi^\# \circ \dbar = d \circ \Phi^\#.
\end{equation*}
Since diffeomorphisms between manifolds preserve the Lipschitz regularity and the vanishing trace condition, $\Phi^\#\dbar\zeta \equiv d(\Phi^\#\zeta)$ is the Euclidean differential of an arbitrary test function in $W^{1,\infty}_0(\Phi^{-1}(\mathscr{D}))$. So, in view of Definition \ref{def: quasiconvex, Euclidean} and \eqref{E quasiconvex, equivalent}, the Riemannian quasiconvexity of $\cc\circ\prz$ on $\mathcal{O}$ is equivalent to the quasiconvexity of 
\begin{equation*}
F(\widehat{x},\widehat{p}) := \sqrt{\det\,g^\E \circ \Phi(\widehat{x})}\, \Phi^\#[\cc\circ\prz] (\widehat{p}),
\end{equation*}
which is defined on $\Phi^{-1}(\mathcal{O}) \subset \R^N$.

Therefore, by Morrey's theory ({\it cf.} \cite{m'}), for $\widehat{p}$ taking the form of a gradient function (identified as the differential via metric duality), namely
\begin{equation*}
\widehat{p}(\widehat{y}) = d^{\bf E}\widehat{\zeta}(\widehat{y}),
\end{equation*} 
the action integral $\int_{\Phi^{-1}(\mathcal{O})}F(\widehat{x},\widehat{p}(\widehat{x}))\,\dd\mathcal{L}^N(\widehat{x})$ is lower semicontinuous with respect to the weak $W^{1,s}$-topology of  the variable $\widehat{\zeta}$. A final change of variables $x \equiv \Phi(\widehat{x})$ yields that
\begin{align*}
\int_{\Phi^{-1}(\mathcal{O})}F\big(\widehat{x},\widehat{p}(\widehat{x})\big)\,\dd\mathcal{L}^N(\widehat{x}) = \int_{\mathcal{O}}\cc\circ\prz \Big( \dbar \big(\Phi_\# \widehat{\zeta}(x)\big) \Big)\,\dmue(x),
\end{align*}
for which $\Phi_\#\widehat{\zeta}$ is again an arbitrary element of $W^{1,\infty}_0(\mathcal{O})$. Since the pushforwards and pullbacks under the diffeomorphism $\Phi$ preserve the weak $W^{1,s}$-topology, we can now conclude the {\em claim} \eqref{claim: lsc}, and hence the theorem follows.     \end{proof}

\section{Discussions}\label{sec: final}

In this work we have proposed, in a primitive fashion, a framework for studying the multi-parameter geometrical variational problems involving differential forms over general fibre bundles. It serves as a first attempt for generalising the programme carried out by Dacorogna--Gangbo in \cite{dg}. Under suitable algebraic and analytic conditions ({\it i.e.}, Definitions \ref{def: horizontal shadow} and \ref{def: quasiconvex, mfd}), we have established the  relaxation and existence theory ({\it cf.} Theorems \ref{thm: equivalence of variational problems} and \ref{thm: existence}) for minimisers of these variational problems. 

Of course, there are many other important aspects of the variational problems lying beyond the scope of this short paper, which shall be left for future investigations. Let us briefly mention three interesting issues here:

\begin{enumerate}
\item
There are other widely studied notions of convexity, other than quasiconvexity, in the calculus of variations over Euclidean domains. These  include rank-one convexity and polyconvexity. Do they admit natural analogues on  Riemannian manifolds, which suitably accounts for the nonlinearity arising from the Riemannian geometry? Also, can we define (Riemannian) quasiconvexity, rank-one convexity and polyconvexity directly for the non-gauged problem $(\wp)$ (see the notations of ``ext-quasiconvexity'', ``ext-rank-one convexity'' and ``ext-polyconvexity'' in \cite{bds})? 
\item
Given the existence theory established in this work, the next natural question is concerned with the regularity  of minimisers. The regularity theory over Euclidean domains has been developed; see {\it e.g.} Evans \cite{e}. Can we develop a (partial) regularity theory over Riemannian manifolds, possibly under the Riemannian quasiconvexity condition?
\item
Can we find physical and engineering models for the application of the relaxation theory developed in this work, for fibres of dimension greater than $1$?
\end{enumerate}

\bigskip
\noindent
{\bf Acknowledgement}.
This work has been done during Siran Li's stay as a CRM--ISM postdoctoral fellow at Centre de Recherches Math\'{e}matiques, Universit\'{e} de Montr\'{e}al and Institut des Sciences Math\'{e}matiques. We  would like to thank these institutions for their hospitality. The author is also grateful to Professors Bob Hardt and Dima Jakobson for stimulating discussions on variational problems over fibre bundles.

\end{document}